\documentclass[11pt,a4paper]{amsart}

\usepackage{fullpage}
\usepackage{amssymb,amsmath,amsthm}
\usepackage{verbatim}
\usepackage{color}
\usepackage
[
  bookmarks = true, 
  pdftitle = \text{K Hughes - A conjecture for arithmetic spherical maximal functions}, 
  colorlinks = false, 
  linkcolor = blue, 
  citecolor = blue
]
{hyperref}

\newcommand{\Z}{\mathbb{Z} }

\newcommand{\N}{\mathbb{N} }

\newcommand{\dimension}{d}

\newcommand{\oddprime}{\mathfrak{p}}
\newcommand{\theprimes}{\mathfrak{P}}

\newtheorem{theorem}{Theorem}
\newtheorem{conjecture}{Conjecture}
\newtheorem{rem}{Remark}

\title{\vspace*{-10mm}A conjecture for arithmetic spherical maximal functions}
\author{K Hughes}
\email{khughes.math@gmail.com}
\date{21 May 2026}

\begin{document}

\begin{abstract}
For 24 years, it has been an open problem to obtain improved bounds, for the maximal function over a sparse sequence of discrete spherical averages, going beyond the range for the full discrete spherical maximal function. 
I formulate a conjecture to characterize the boundedness of such maximal functions and state a theorem in support of it. 
\end{abstract}

\maketitle



We begin with some notation so that we may precisely formulate the problem. 
Let \(d \geq 4\). 
Define the arithmetic spherical averages
\[
A_\lambda f(x) = \frac{1}{\#\{y \in \Z^\dimension: |y|^2=\lambda\}} \sum_{y \in \Z^\dimension: |y|^2=\lambda} f(x-y)
\]
for $\lambda \in \N$, \(f \in \ell^1(\Z^d)\) and \(x \in \Z^d\).  
For $\Lambda \subset \N$, define the arithmetic spherical maximal function 
\[
M_\Lambda f(x) := \sup_{\lambda \in \Lambda} |A_\lambda f(x)|. 
\]
In \cite{MSW:spherical}, it was proved that \(M_{\N}\) is bounded on \(\ell^p(\Z^d)\) if and only if \(d \geq 5\) and \(p>\frac{d}{d-2}\). The restricted weak-type endpoint result was given in \cite{Ionescu}. Let \(\Lambda \subseteq \N\). Since \(M_{\N}f(x) \geq M_{\Lambda}f(x)\) for all \(f \in \ell^p(\Z^d)\) and \(x \in \Z^d\), \cite{MSW:spherical} implies that \(M_{\Lambda}\) is bounded on \(\ell^p(\Z^d)\) when \(d \geq 5\) and \(p>\frac{d}{d-2}\). 
A natural question is: \textsf{\textcolor{blue}{Is this range of \(\ell^p(\Z^d)\) spaces sharp for a given sequence \(\Lambda\)?}}

By analogy with the continuous Euclidean lacunary spherical maximal function \cite{Calderon,CW}, the folklore conjecture was that, for each lacunary sequence $\Lambda \subset \N$, $M_\Lambda$ is bounded on $\ell^p(\Z^d)$ for all $p>1$ and $d \geq 5$. I also conjectured this is true when \(d=4\) for sequences avoiding the 2-adic obstruction when \(d=4\). 
The first progress towards this question was made in \cite{Hughes:sparse}. Subsequently, \cite{Cook:sparse, Cook:Birch} showed that there are infinite sequences that are bounded on \(\ell^p(\Z^d)\) for \(d \geq 5\) and all \(p>1\). At my behest, \cite{KLM} combined my ideas and methods introduced in \cite{Hughes:sparse} and \cite{Hughes:restricted} to prove that -- for \emph{any} lacunary sequence \(\Lambda \subset \N\) -- \(M_{\Lambda}\) is bounded on \(\ell^p(\Z^d)\) when \(p>\frac{d-2}{d-3}\) and \(d \geq 5\). 
Such maximal functions have been studied further in \cite{CH, AM, MSW:dimfree, LSW, Zhang, BCSS}. 

While writing up \cite{Hughes:sparse}, Prof. Jacek Zienkiewicz showed me a probabilistic construction proving the existence of lacunary sequences \(\Lambda \subset \N\) whose associated arithmetic spherical maximal function \(M_{\Lambda}\) is unbounded near $\ell^1(\Z^\dimension)$ for each \(d \geq 5\). 
His construction disproved the folklore conjecture. 
I quickly adapted the construction to show that, for each \(p<\frac{d}{d-1}\) and \(d \geq 4\), there exist arbitrarily fast growing sequences, say \(\Lambda\), in \(\N\) such that \(M_{\Lambda}\) is unbounded on \(\ell^p(\Z^d)\). This showed that intuition from the continuous setting did not apply; see \cite{CH}. 
Yet, the intuition from the continuous setting can apply and the folklore conjecture can be true for averages defined by objects similar to spheres, but with substantially less arithmetic, such as the equation \(\lambda = \sum_{i=1}^d \lfloor x_i^c \rfloor\) where \(c\) is not integral; see \cite{ILMS}. 

Given the mixture of positive and negative results, the above question becomes: \textsf{\textcolor{blue}{For a given sequence \(\Lambda \subset \N\), can we characterize the range of \(\ell^p(\Z^d)\) spaces on which \(M_{\Lambda}\) is bounded?}}
The point of this note is to conjecture a characterization. 

To set this up: Let \(\theprimes\) denote the set of rational primes. 
For each prime \(\oddprime \in \theprimes\), let \(\delta_\oddprime(\Lambda)\) denote the \(\oddprime\)-adic upper Minkowski dimension of \(\Lambda\). 
Let \(\delta_\infty(\Lambda)\) denote the smallest parameter \(\delta \geq 0\) such that there exists a constant \(C_\delta\) -- independent of \(N \in \N\) -- with
\(
\#(\Lambda \cap [N,2N]) 
\leq C_\delta 
N^{\delta} 
\) uniformly for all \(N \in \N\).  
These dimensional quantities are analogous to those in \cite{SWW}. 
Define the exponent 

\[
\eta(\Lambda)
:= 
\max\Big\{ \sup_{\oddprime \in \theprimes} \{1+\frac{\delta_\oddprime(\Lambda)}{\dimension-1}\} ,  1+\frac{2\delta_\infty(\Lambda)}{\dimension-2} \Big\}
.\]

We have the following unboundedness result (from 2015-6) whose proof I will defer to another paper (though its statement makes the proof obvious). 
\begin{theorem}\label{thm:main}
Let \(\Lambda \subset \N\). 
\( M_{\Lambda} \) is unbounded on \(\ell^p(\Z^\dimension)\) for  all \(p<\eta(\Lambda)\). 
\end{theorem}

Finally, I conjecture the following. 
\begin{conjecture}
\label{conj:2016}
Let \(\Lambda \subset \N\). 
Let \(d \geq 5\). 
\( M_{\Lambda} \) is bounded on \(\ell^p(\Z^\dimension)\) for all \(p>\eta(\Lambda)\). 
\end{conjecture}

\begin{rem}
In particular, Conjecture~\ref{conj:2016} states that any lacunary arithmetic spherical maximal function would be bounded on \(\ell^p(\Z^\dimension)\) for \( p > \frac{\dimension}{\dimension-1} \) and \(d \geq 5\). 
\end{rem}

\begin{rem}
Theorem~\ref{thm:main} and Conjecture~\ref{conj:2016} generalize to non-degenerate positive integral quadratic forms and to higher degree homogeneous integral forms. 
These generalizations will be investigated in a subsequent paper. 
\end{rem}

\bibliographystyle{amsalpha}

\end{document}